\documentclass[11pt]{article}
\usepackage{amsthm,amssymb,amscd,amsmath,latexsym,indentfirst,color,amsfonts}
\usepackage[mathcal]{eucal}
\theoremstyle{plain}
\newtheorem{theorem}{Theorem}[section]
\newtheorem{prop}[theorem]{Proposition}
\newtheorem{coro}[theorem]{Corollary}
\newtheorem{lem}[theorem]{Lemma}
\newtheorem{defi}[theorem]{Definition}
\theoremstyle{remark}
\newtheorem{rem}[theorem]{Remark}

\DeclareMathOperator\Lie{Lie}
\DeclareMathOperator\PGL{PGL}\DeclareMathOperator\GL{GL}\DeclareMathOperator\Sp{Sp}\DeclareMathOperator\SL{SL}
\DeclareMathOperator\Ad{Ad} \DeclareMathOperator\Aut{Aut}
\DeclareMathOperator\Hom{Hom}\DeclareMathOperator\Int{Int}
\DeclareMathOperator\ord{ord}\DeclareMathOperator\ad{ad}
\DeclareMathOperator\Ker{Ker}\DeclareMathOperator\an{and}

\def\bt{\beta} \def\dt{\delta}
\def\al{\alpha}
\def\rt{\rightarrow}\def\st{\subset}\def\ddt{\Delta}

\def\ot{\otimes}\def\op{\oplus}\def\si{\sigma}\def\Ga{\Gamma}\def\wh{\widehat}\def\ga{\gamma}\def\wt{\widetilde}
\def\mb{\mathbb}\def\va{\varphi}

\newcommand{\be}{\begin {equation}}
\newcommand{\ee}{\end {equation}}
\newcommand{\bp}{\begin {proof}}
\newcommand{\ep}{\end {proof}}
\newcommand{\bee}{\begin {equation*}}
\newcommand{\eee}{\end {equation*}}
\newcommand{\lb}{\label}

\begin{document}
\title{The Weyl group of the fine grading of $sl(n,\mathbb{C})$ associated with tensor product of generalized Pauli matrices \footnote{Research supported by NSFC Grant No.10801116 and by
'the Fundamental Research Funds for the Central Universities'}}
\author{Gang HAN \\
{Department of Mathematics, Zhejiang University,}\\{Hangzhou 310027, China}\\
{E-mail: mathhg@hotmail.com}\\Telephone: 086-0571-87953843\\Fax:
086-0571-87953794 }
\date{January 24, 2011}
\maketitle

{\small \noindent \textbf{Abstract}.  We consider the fine grading
of $sl(n,\mb C)$ induced by tensor product of generalized Pauli
matrices in the paper. Based on the classification of maximal
diagonalizable subgroups of $PGL(n,\mb C)$ by Havlicek, Patera and
Pelantova, we prove that any finite maximal diagonalizable subgroup
$K$ of $PGL(n,\mb C)$ is a symplectic abelian group and its Weyl
group, which describes the symmetry of the fine grading induced by
the action of $K$, is
 just the isometry group of the symplectic abelian group $K$.
For a finite symplectic abelian group, it is also proved that its isometry group is always generated by the transvections
contained in it.


 }

\section{Introduction}
 \setcounter{equation}{0}\setcounter{theorem}{0}
The study of gradings of Lie algebras and the symmetries of those
gradings is an active research area in recent decades, which are
interesting to both mathematicians and physicians. In physics, Lie
algebras usually play the role as the algebra of infinitesimal
symmetries of a physical system. Knowledge about the gradings of a
Lie algebra will greatly help us to understand better the structure
of the Lie algebra. Study of the symmetries of those gradings offers
a very important tool for describing symmetries in the system of
nonlinear equations connected with contraction of a Lie algebra (see
e.g. \cite{jpt}).

 Besides the
famous Cartan decomposition for semisimple Lie algebras, another
well-known example of grading is the grading of $sl(n,\mb C)$ by the
adjoint action of the Pauli group $\Pi_n$ generated by the $n\times
n$ generalized Pauli matrices, which decomposes $sl(n,\mb C)$ into
direct sum of $n^2-1$ one-dimensional subspaces, each of which
consists of semisimple elements.

 Let $L$ be a complex simple Lie algebra. Let $\Aut(\L)$ and $\Int(L)$ be respectively the automorphism group and inner automorphism group of $L$, which are
both algebraic groups. A subgroup of $\Aut(\L)$ or $\Int(L)$ is
called \textit{diagonalizable} if it is abelian and consists of
semisimple elements. It is not hard to see that there is a natural
1-1 correspondence between gradings of $L$ and diagonalizable
subgroups of $\Aut(\L)$ (see Section 4). A grading is called
\textit{inner} if the respective diagonalizable subgroup is in
 $\Int(\L)$. A grading (resp. inner grading) of $L$ is called fine if it could not be further refined by any other grading (resp. inner grading).
 Among the gradings of a Lie algebra, fine (inner) gradings are especially
important. It was shown in \cite{pz} that the fine gradings of
simple Lie algebras correspond to maximal diagonalizable subgroups
(which were called MAD-groups in \cite{pz}) of $\Aut(\L)$. Then fine
inner grading of simple Lie algebras corresponds to maximal
diagonalizable subgroups of $\Int(\L)$. Given a fine (inner) grading
$\Ga$, one can define naturally its Weyl group (see Definition 2.3
of \cite{hg}) to describe its symmetry. Assume $K$ is the maximal
diagonalizable subgroup corresponding to $\Ga$, then one can show
 that its Weyl group is isomorphic to the Weyl group of $K$. See Proposition 2.4 and Corollary 2.6 of \cite{hg}.

After many mathematicians and physicians' contribution,  the
classification of fine gradings of all the simple Lie algebras are
almost done. For example, it can be found in \cite{e} the
classification of fine gradings of all the classical simple Lie
algebras over an algebraically closed field of characteristic 0.
People have also known a lot about the fine gradings for exceptional
simple Lie algebras, see \cite{k} for a survey of such results. For
the Weyl group of a fine inner grading of a simple Lie algebra $L$,
if the grading is Cartan decomposition (in which case the
corresponding maximal diagonalizable subgroup $K$ of $\Int(L)$ is
just the maximal torus), then it is well-known that the Weyl group
is a finite group generated by reflections; in other cases there is
no general result by far. The next step is to study the case $K$ is
discrete, and people have made some explorations in the case
$L=sl(n,\mb C)$ .

Recall that $PGL(n,\mb C)$ is the inner automorphism group of
$sl(n,\mb C)$. Let us first review the classification of maximal
diagonalizable subgroups of $PGL(n,\mb{C})$, which correspond to
fine inner gradings of $sl(n,\mb C)$. Let $\Pi_n$ be the Pauli group
of $GL(n,\mb{C})$ and $D_n$ be the subgroup of diagonal matrices of
$GL(n,\mb{C})$. Let
 $\texttt{P}_n$  and $\texttt{D}_n$ be the respective images of $\Pi_n$ and $D_n$ in $PGL(n,\mb{C})$ under the adjoint
action on $M(n,\mb{C})$. Assume $n=k l_1\cdots l_t$ and each $l_i$
divides $l_{i-1}$. The group $D_k\ot\Pi_{l_1}\ot \cdots
\ot\Pi_{l_t}$ consists of all those elements $A_0\ot A_1\ot\cdots\ot
A_t$ with $A_0\in D_k$ and $A_i\in\Pi_{l_i}$ for $1\le i\le t$. The adjoint action of
$D_k\ot\Pi_{l_1}\ot \cdots \ot\Pi_{l_t}$ on
$$M(k,\mb{C})\ot M(l_1,\mb{C})\ot\cdots\ot M(l_t,\mb{C})\cong M(n,\mb{C})$$ induces the embedding
$$ \texttt{D}_k\times \texttt{P}_{l_1}\times \cdots \times \texttt{P}_{l_t}\hookrightarrow PGL(n,\mb{C}). $$ If we identify
$ \texttt{D}_k\times \texttt{P}_{l_1}\times \cdots \times
\texttt{P}_{l_t}$ with its image, then it was shown by Havlicek,
Patera and Pelantova in Theorem 3.2 of \cite{hpp} that any maximal
diagonalizable subgroup $K$ of $PGL(n,\mb{C})$ is conjugate
 to one and only one of the $ \texttt{D}_k\times \texttt{P}_{l_1}\times \cdots \times \texttt{P}_{l_t}$.

Let $K$ be a discrete maximal diagonalizable subgroup of $PGL(n,\mb
C)$. Then $K\cong\texttt{P}_{l_1}\times\texttt{P}_{l_2}\times \cdots
\times \texttt{P}_{l_t}$ where $n=l_1 l_2\cdots l_t$ and each $l_i$
divides $l_{i-1}$. Then the fine grading induced by $K$ also
decomposes $sl(n,\mb C)$ into $n^2-1$ one-dimensional subspaces,
each of which consists of semisimple elements. We will show in
Section 5 that there
 is a nonsingular anti-symmetric pairing $<,>$ on $K$,
 such that
 $(K,<>)$ is a nonsingular symplectic abelian group (see Definition \ref{d8}). Moreover the pairing $<,>$ is invariant under the Weyl group of
 $K$. It is shown in Proposition \ref{d6} that there is a one-to-one correspondence between conjugacy classes
of finite maximal diagonalizable subgroups of $\PGL(n,\mathbb{C})$
and finite symplectic abelian groups of order $n^2$. The following
theorem about the structure of the isometry group of a finite
nonsingular symplectic abelian group is Theorem \ref{c3}. For the
definition of a transvection on a symplectic abelian group, see
Definition \ref{d4}.

\begin{theorem}\lb{cc3}
Let $(H,<,>)$ be a finite nonsingular symplectic abelian group. Then
its isometry group is generated by the set of transvections in it.
\end{theorem}

 If $K=P_n$, then in the important paper \cite{jpt} the authors showed that the respective Weyl group is $\SL(2,\mb Z_n)$. If $m=p^2$ with $p$ a
prime and $K=\texttt{P}_p\times \texttt{P}_p$, then in \cite{pst}
the authors proved that the respective Weyl group is
$\Sp(4,\mathbb{Z}_p)$. Next in \cite{hg} we dealt with the case
$n=m^k$ ($m$ may not be a prime) and $K=\texttt{P}_m^{k}$ is the
$k$-fold direct product of $\texttt{P}_m$, and proved that the Weyl
group is isomorphic to $\Sp(2k,\mathbb{Z}_n)$ and is generated by
transvections. Then, in this paper we deal with the general case
that $K$ is an arbitrary discrete maximal diagonalizable subgroup,
and prove the following result in Theorem \ref{d7} generalizing the
previous result.

\begin{theorem}\lb{dd9}

Let $K$ be a finite maximal diagonalizable subgroup of
$G=\PGL(n,\mathbb{C})$ and $W_G(K)$ be its Weyl group. Then $W_G(K)$
equals the isometry group of $(K,<,>)$, and is generated by the set
of transvections in it.
\end{theorem}

 The paper is organized as follows. The definition and classification of
finite nonsingular symplectic abelian groups will be reviewed in
Section 2. Then in Section 3 we will define transvections on a
finite symplectic abelian group and prove Theorem \ref{cc3}. Next in
section 4 we will review the definitions of the grading of a simple
Lie algebra and prove some important properties. Then in Section 5
we will define the anti-symmetric pairing on any finite maximal
diagonalizable subgroup of $\PGL(n,\mathbb{C})$ and prove the $1-1$
correspondence between conjugacy classes of finite maximal
diagonalizable subgroups of $\PGL(n,\mb C)$ and nonsingular
symplectic abelian groups of order $n^2$. In the last section,
Theorem \ref{dd9} will be proved.\bigskip

Finally we introduce some notations in the paper.

For a finite set $S$, we will use $|S|$ to denote its cardinality.

For any $n\in \mb Z_+$, let $\mb Z_n=\mb Z/n\mb
Z=\{\bar{0},\bar{1},\cdots,\overline{n-1}\}$. For simplicity
 we will just use $i$ to denote $\bar{i}$ for $i=0,1,\cdots,n-1$.

Let $\omega_n=e^{2\pi i/n}$ and
$C_n=\{\omega_n^i|i=0,1,\cdots,n-1\}$ be the cyclic group of order
$n$ generated by $\omega_n$. Sometimes we will identify $\mb Z_n$
with $C_n$ by mapping $i$ to $\omega_n^i$.
\bigskip

\centerline{\textit{Acknowledgments}}

The research is supported by NSFC Grant No.10801116 and by 'the
Fundamental Research Funds for the Central Universities'. It is
finished during the author's visit at the department of mathematics
in MIT in 2011. He acknowledges the hospitality of MIT, and would
like to take this opportunity to heartily thank David Vogan for
drawing his attention to this problem and for Vogan's great
generosity in sharing his immense knowledge with him during the
research. Proposition \ref{d0} is due to him.

\section{Classification of finite nonsingular symplectic abelian groups
}
 \setcounter{equation}{0}\setcounter{theorem}{0}
We will follow the definition of symplectic abelian groups in
\cite{ka}, which is defined with respect to any field. But for our
purpose we will always assume the field to be $\mb C$,
and we will write our abelian groups additively in Section 2 and 3.

Let $H$ be an abelian group. Recall that an abelian group is
automatically a $\mb Z$-module.
\begin{defi}\lb{d8}
 A map  $$<,>:H\times H\rt \mb C^{\times}$$ is called a pairing of $H$ into $\mb C^{\times}$ if
$<,>$ is $\mb Z$-bilinear. The pairing is called anti-symmetric if
for all $a,b\in H$, $$<a,b>=<b,a>^{-1}.$$ An anti-symmetric pairing
$<,>$ is called nonsingular if $<a,b>=1$ for any $b\in H$ implying
$a=0$.
\end{defi}
\begin{defi}
Assume $<,>$ is an anti-symmetric pairing of $H$ into $\mb C^{\times}$. Then $(H,<,>)$
is called a \textit{symplectic abelian group}. A symplectic abelian group $(H,<,>)$ is said to be
\textit{nonsingular} if $<,>$ is nonsingular.
\end{defi}

Now assume that $(H,<,>)$ is a nonsingular symplectic abelian group.

 A subgroup $H_0$ of $H$ is called a \textit{nonsingular symplectic
 abelian subgroup} if the restriction $<,>|H_0$ is nonsingular.

Two subgroups $H_1$ and $H_2$ of $H$ are said to be
\textit{orthogonal}, written $H_1\perp H_2$, if $<a,b>=1$ for any
$a\in H_1, b\in H_2$.

 Two symplectic abelian
groups are said to be \textit{isometric} if there is a group isomorphism
between them preserving the respective pairings.

 If
$H_1,H_2,\cdots,H_n$ is a family of nonsingular symplectic
 abelian subgroups of $H$ such that
$$H=H_1\oplus H_2\oplus\cdots\oplus H_n$$ and $$H_i\perp H_j,\ i\neq j,$$
then we will say that $H$ is the orthogonal direct sum of symplectic
abelian subgroups $H_1,H_2,\cdots,H_n$.

Assume $n\in \mb Z^{+}$ and $n>1$. If a pair of elements $u,v\in H$
of order $n$ satisfying $<u,v>=\omega_n$, then we call $(u,v)$ a
\textit{hyperbolic pair} of order $n$ in $H$. Let \be \mb H_n=\mb
Z_n\times \mb Z_n\lb{d3}\ee and $<(i,j),(k,l)>=\omega_n^{il-jk}$ be
the pairing on $\mb H_n$, which is clearly nonsingular and
anti-symmetric. Then $(\mb H_n,<,>)$ (or just $\mb H_n$) is a
nonsingular symplectic abelian group, called the \textit{hyperbolic
group} of rank $n$. Note that in \cite{ka} the rank $n$
for a hyperbolic group is assumed to be a power of a prime, but we
will not have this restriction in this paper. Let $u_1=(1,0)\in \mb
H_n,v_1=(0,1)\in\mb H_n$. Then the hyperbolic pair $(u_1,v_1)\in\mb
H_n^2$ is called the \textit{standard hyperbolic pair} of $\mb H_n$.

\begin{lem}\lb{e}
Let $H$ be a finite symplectic abelian group. If $a,b\in H$ are both
of order $n$ and $<a,b>=\omega_n$, then $a$ and $b$ generate a
subgroup $K$ isometric to $\mb H_n$.
\end{lem}
\bp Assume for some $i,j\in \mb Z$, $ia+jb=0$. Then
$<a,ia+jb>=\omega_n^j=1$ thus $n|j$. Similarly $n|i$. So
$K=\{ia+jb|i,j\in \mb Z_n\}$. Then $$K\rt \mb H_n,\ ia+jb\mapsto
(i,j)$$ is clearly an isometry of symplectic abelian groups. \ep
\begin{lem}\lb{c}
If $m$ and $n$ are relatively prime, then $\mb H_{mn}\cong \mb H_m
\op \mb H_n$ as symplectic abelian groups.
\end{lem}
\bp Let $(u_1,v_1)\in \mb H_m^2$ (resp. $(u_2,v_2)\in \mb H_n^2$) be
the standard hyperbolic pair for $\mb H_m$ (resp. $\mb H_n$). Let
$$a=u_1+u_2\in \mb H_m \op \mb H_n,\ b=v_1+v_2\in \mb H_m
\op \mb H_n.$$ The order of $a$ and $b$ are both $mn$ as $m$ and $n$
are relatively prime.

One has $<a,b>=<u_1,v_1><u_2,v_2>=\omega_m
\omega_n=\omega_{mn}^{m+n}$. As $mn$ and $m+n$ are also relatively
prime, $<a,ib>=\omega_{mn}$ for some integer $i$. Clearly
$\ord(ib)$, the order of $ib$, is still $mn$. Then by Lemma \ref{e},
the subgroup of $\mb H_m \op \mb H_n$ generated by $a$ and $ib$ is
isometric to $\mb H_{mn}$. As $|\mb H_{mn}|=|\mb H_m \op \mb H_n|$,
$$\mb H_m \op \mb H_n\rt \mb H_{mn},\ ja+k(ib)\mapsto (j,k)$$ is an
isometry of symplectic abelian groups. \ep
\begin{lem}\lb{x}
Let $(H,<,>)$ be a symplectic abelian group. Let $a,b\in H$ with
$ord(a)=i,\ ord(b)=j$. Assume $<a,b>=x$.

(1) If $l$ is the minimal positive integer such that $x^l=1$, then
$l|i$ and $l|j$.

(2) If $i,j$ are relatively prime then $x=1$.
\end{lem}
\bp (1) As $x^i=<ia,b>=<0,b>=1$, one has $l|i$. Similarly $l|j$.

(2) Apply (1). \ep

For any prime $p$ dividing $|H|$, we will always denote the $p$-Sylow subgroup of $H$ by $H(p)$.
\begin{theorem}\lb{d1}[Lemma 1.6 and Theorem 1.8 of \cite{ka}]
Let $(H,<,>)$ be a finite nonsingular symplectic abelian group. Then

(1) $H$ is an orthogonal direct sum of all its Sylow subgroups.

(2) Assume that $H(p)$ is a $p$-Sylow subgroup of $H$. Then $$H(p)\cong \mb H_{p^{r_1}}\op \mb H_{p^{r_2}}\op\cdots \op \mb
H_{p^{r_s}}$$ for some positive integers $r_1,r_2,\cdots,r_s$ with
$r_i\geq r_{i+1}$.

(3) $H$ is an orthogonal direct sum of hyperbolic subgroups $\mb
H_n$, where each $n$ is a power of some prime.
\end{theorem}

\bp (1) It is proved in Lemma 1.6 of \cite{ka}. In fact it is an
easy consequence of Lemma \ref{x} (2).

(2) This is proved in Theorem 1.8 of
\cite{ka} implicitly.

(3)  It is a consequence of (1) and (2).
It is also proved in Theorem 1.8 of
\cite{ka}.
\ep

\begin{coro}\lb{d9}
Let $(H,<,>)$ be a nonidentity finite nonsingular symplectic abelian
group. Then there exists positive integers $l_1,l_2,\cdots,l_k$ with
each $l_i|l_{i-1}$ such that \be H\cong \mb H_{l_1}\op \mb
H_{l_2}\op \cdots \op \mb H_{l_k}\lb{b}.\ee Such positive integers
are uniquely determined by $H$.
\end{coro}
\bp The existence of such positive integers and the isometry
(\ref{b}) follow directly from (1) and (2) of Theorem \ref{d1}.
According to the structure theorem of finite abelian groups, such
positive integers are uniquely determined by $H$. \ep

 Let $(H,<,>)$ be a nonsingular symplectic abelian group and $H_0$ be a finite nonsingular symplectic subgroup. For any $a\in H$ of
 order $n$ and any $j\in \mb
 Z$, define
 \be \omega_n^j\cdot a=_{def}ja.
 \lb{t1}\ee
  Assume that $H_0\cong \mb H_{l_1}\op \mb H_{l_2}\op \cdots \op
\mb H_{l_k}$ with $l_i|l_{i-1}$ for all $i$. For $i=1,\cdots,k$ let
$(a_i,b_i)$ be a hyperbolic pair of order $l_i$ for $\mb H_{l_i}$.
Let \be\pi:H\rt H_0,\ c\mapsto \sum_{i=1}^k(<c,b_i>\cdot
a_i-<c,a_i>\cdot b_i).\lb{p}\ee Let $H_0^{\perp}=\{a\in H|<a,b>=0,
\forall\ b\in H_0\}$.
\begin{prop}
(1)$H_0^{\perp}$ is a symplectic subgroup of $H$ and $H=H_0\op
H_0^{\perp}$ is a direct sum of nonsingular symplectic abelian
subgroups.

(2)The map $\pi$ is independent of the hyperbolic pairs $(a_i,b_i)$
chosen.
\end{prop}
\bp
 (1) First it is clear that $H_0^{\perp}$ is a subgroup of $H$.
$H_0\cap H_0^{\perp}=0$ as \\ $<,>|H_0$ is nonsingular. Let $c\in
H$. Assume $<c,a_i>=\omega_{l_i}^{t_i}$ for $i=1,\cdots,k$, then
\bee\begin{split}<\pi(c),a_i>&=<-<c,a_i>b_i,a_i>\\&=<-t_i
b_i,a_i>\\&=<b_i,a_i>^{-t_i}=\omega_{l_i}^{t_i}.\end{split}\eee So
$$<c,a_i>=<\pi(c),a_i>$$ for
$i=1,\cdots,k$. Similarly $$<c,b_i>=<\pi(c),b_i>$$ for
$i=1,\cdots,k$. Thus one has $c-\pi(c)\in H_0^{\perp}$ and \be
c=\pi(c)+(c-\pi(c))\in H_0+H_0^{\perp}.\lb{q}\ee So $ H=H_0\op
H_0^{\perp}.$ It is clear that $<,>|H_0^{\perp}$ must also be
nonsingular as $<,>$ is nonsingular. Thus $H=H_0\op H_0^{\perp}$ is a direct sum of
nonsingular symplectic abelian subgroups.

 (2) Let $\pi^{'}:H\rt
H_0$ be defined as in (\ref{p}) with respect to another choice of
hyperbolic pairs $(a_i^{'},b_i^{'})$ for each $\mb H_{l_i}$. Then
for any $c\in H$ one also has $$c=\pi^{'}(c)+(c-\pi^{'}(c))\in
H_0\op H_0^{\perp}.$$ Comparing to (\ref{q}), as $H=H_0\op
H_0^{\perp}$ is a direct sum, one must have \\ $\pi(c)=\pi^{'}(c)$
for any $c\in H$. \ep We call the map $\pi:H\rt H_0$ defined in
(\ref{p}) the \textit{projection} of $H$ onto $H_0$.
\section{Transvections and isometry groups of finite nonsingular symplectic abelian groups
}
 \setcounter{equation}{0}\setcounter{theorem}{0}
 Let $(H,<,>)$ be a finite nonsingular symplectic abelian group.

\begin{defi}
 Let $\Sp(H)$ be the set of isometries of $H$ onto itself. Then $\Sp(H)$ is clearly a group, called the isometry group
 of $H$.
\end{defi}

If $H=H_1\op H_2$ is a direct sum of nonsingular symplectic abelian
subgroups then it is clear that
$$\Sp(H_1)\times \Sp(H_2)\rt \Sp(H),\
(\phi,\nu)(a,b)=(\phi(a),\nu(b))$$ embeds $\Sp(H_1)\times \Sp(H_2)$
as a subgroup of $\Sp(H)$. For any symplectic subgroup $H_0$ of $H$,
as $H=H_0\op H_0^{\perp}$, we will always regard $\Sp(H_0)$ as a
subgroup of $\Sp(H)$ by the embedding
$$\Sp(H_0)\hookrightarrow \Sp(H_0)\times \Sp(H_0^{\perp})\st
\Sp(H),\ \phi\mapsto (\phi,1),$$ where 1 denotes the identity map on
$H_0^{\perp}$.

\begin{prop}
 Let $H$ be a finite nonsingular symplectic abelian group. Then $\Sp(H)$ acts transitively on the set of hyperbolic pairs in $H$ with the same order.
\end{prop}
\bp Assume $(a_1,b_1)$ and  $(a_2,b_2)$ are two hyperbolic pairs in
$H$ with order $n$. Let $H_i=Span (a_i,b_i)$ for $i=1,2$. Then $H_i$
is a symplectic subgroup of $H$ isometric to $\mb H_n$. Let
$\phi_1:H_1\rt H_2$ be the isometry such that $\phi_1(a_1)=a_2,\
\phi_1(b_1)=b_2.$

One has $H=H_1\op H_1^{\perp}$ and $H=H_2\op H_2^{\perp}$. By
Theorem \ref{d1} (3), $H_1^{\perp}\cong H_2^{\perp}$. Fix an
isometry $\phi_2:H_1^{\perp}\rt H_2^{\perp}$. Then
$$\phi=(\phi_1,\phi_2):H_1\op H_1^{\perp}\rt H_2\op H_2^{\perp},\
(c,d)\mapsto (\phi_1(c),\phi_2(d))$$ is in $\Sp(H)$ and maps
$(a_1,b_1)$ to  $(a_2,b_2)$. Thus $\Sp(H)$ acts transitively on the
set of hyperbolic pairs in $H$ with the same order. \ep

  Let $$\wh H=_{def}Hom(H,\mb C^{\times})$$ be the abelian group of characters of $H$.

For any $a\in H$, define
 $$\ga_a:H\rt \mb C^{\times},\ b\mapsto <a,b>.$$

 \begin{lem} The map \be \varphi:H\rt \wh H,\ a\mapsto \ga_a\lb{g}\ee is an isomorphism of abelian groups.\end{lem}
 \bp
As $<,>$ is $\mb Z$-bilinear, $\va$ is $\mb Z$-linear. As $<,>$ is
nonsingular, $\va$ is one-to-one, thus is onto as $|H|=|\wh H|$.
 \ep
  \begin{defi}\lb{a5}
   For any $\ga\in \wh H$, let $$\ga^*=_{def}\varphi^{-1}(\ga).$$
  \end{defi}
Then for any $a\in H$, $\ga(a)=<\ga^*,a>$.
  As (\ref{g}) is an isomorphism, $\ga$ and $\ga^*$ have the same
  order, and for $\ga_1,\ga_2\in \wh H$,
  \be (\ga_1+\ga_2)^*=\ga_1^*+\ga_2^*.\lb{b3}\ee

Recall that for a vector space $V$ over $\mb C$ (or other field), a linear map
$\phi\in \GL(V)$ is called a transvection if
$$\phi(v)=v+\lambda(v)u$$ for some $\lambda\in V^*, u\in V$ satisfying
$\lambda(u)=0$. If $V$ is a symplectic vector space with $<,>$ the
anti-symmetric pairing on it, then any transvection of $V$
preserving the form $<,>$ must be of the form
\be\phi(v)=v-k<u,v>u\lb{d5}\ee for some $u\in V,\ k\in \mb C$. One
knows that $\SL(V)$ and $\Sp(V)$ are both generated by their
transvections. One can define transvections for a symplectic abelian
group analogously.

For any $b\in H$ with $b\neq 0$, assume $\ord(b)=m$. For any $a\in
H$, $<b,a>=\ga_b(a)$ takes value in the cyclic group
$C_m=\{\omega_m^i|i=0,1,\cdots,m-1\}$. Recall the convention in
(\ref{t1}).
\begin{defi}\lb{d4} For any $b\in H$ with $b\neq 0$ and $k\in\mb Z$,
define a homomorphism \be s_{b,k}:H\rt H,\ a\mapsto a-k(<b,a>\cdot
b),\lb{g3}\ee and call it a transvection on $H$. Using the
identification $\varphi$ of $H$ and $\wh H$, for any $\ga\in \wh H$
define $ s_{\ga,k}= s_{\ga^*,k}$. Denote $s_{b,1}$ (resp.
$s_{\ga,1}$) by $s_{b}$ (resp. $s_{\ga}$) for simplicity.
 \end{defi}
 Then one
has \be s_\ga(a)=a-\ga(a)\ga^*.\lb{b9}\ee

\begin{lem} (1) For any $b\in H$ with $b\neq 0$ and any $k,j\in\mb Z$, one has
 \be  s_{b,0}=1,\lb{g0}\ee \be s_{b,k}
 s_{b,j}= s_{b,k+j}, \lb{g1}\ee and
\be s_{b,k}^{-1}= s_{b,-k}\lb{g2}\ee

(2) If $\ord(b)=m$, them $\{ s_{b,k}|k\in\mb Z\}$ is a cyclic
  group of order $m$ generated by $ s_{b,1}$.
\end{lem}

\bp (1) (\ref{g0}) follows from the definition (\ref{g3}). For any
$a\in H$, \bee\begin{split} s_{b,k}
 s_{b,j}(a)&= s_{b,k}(a-j(<b,a>\cdot b))\\&=(a-j(<b,a>\cdot
b))-k(<b,a-j(<b,a>\cdot b)>\cdot b)\\&=a-(k+j)(<b,a>\cdot b)\\&=
s_{b,k+j}(a),\end{split}\eee So (\ref{g1}) holds. Then (\ref{g2})
follows from (\ref{g0}) and (\ref{g1}).

(2) It follows from (1). \ep

\begin{lem}
 One has $ s_{b,k}\in \Sp(H)$, where $b\in H$ with $b\neq 0$ and $k\in\mb
 Z$.
\end{lem}
\bp By (2) of last lemma one only need to show $ s_b\in \Sp(H)$. By
(\ref{g3}) and (\ref{g2}) it is clear that $ s_b$ is a $\mb
Z$-linear isomorphism of $H$, so we only need to prove that $ s_b$
preserves $<,>$. Assume $a,c\in H$ and $<b,a>=\omega_m^i$,
$<b,c>=\omega_m^j$. Then \bee
\begin{split}
< s_b(a), s_b(c)>&=<a-<b,a>b,c-<b,c>b>\\
&=<a-ib,c-jb>\\&=<a,c> <b,c>^{-i}<a,b>^{-j}\\
&=<a,c> \omega_m^{j(-i)}\omega_m^{(-i)(-j)}\\
&=<a,c>.
\end{split}
\eee

\ep

Let \be Q(H)=_{def}< s_{b,k}|0\neq b\in H,k\in \mb Z>\ee be the
subgroup of $\Sp(H)$ generated by all the transvections. It is clear
that $Q(H)$ is generated by those $ s_b$ with $b\in H$ and $b\neq
0$. For any element $b\neq 0$ in a nonsingular symplectic subgroup
$H_0$ of $H$, $ s_b\in Q(H_0)$ can be regarded as in $Q(H)$ since
$b\in H$. Thus $Q(H_0)$ can be naturally regarded as a subgroup of
$Q(H)$.

Let $\GL(2,\mathbb{Z}_n)$ be the group of $2\times 2$ invertible
matrices in $M(2,\mathbb{Z}_n)$. Let $$\SL(2,\mb Z_n)=\{A\in
\GL(2,\mathbb{Z}_n)|\det(A)=1\in\mb Z_n\}.$$ Let $$J=\left(
             \begin{array}{cc}
               0 & 1\\
               -1 & 0 \\
             \end{array}
           \right)\in M(2,\mb Z_n)$$

and $$\Sp(2,\mb Z_n)=\{A\in \GL(2,\mathbb{Z}_n)|A^t J A=J\}.$$ It is
easily verified that $\SL(2,\mb Z_n)=\Sp(2,\mb Z_n)$.

\begin{lem}\lb{h}
One has $\Sp(\mb H_n)\cong \Sp(2,\mb Z_n)=\SL(2,\mb Z_n)$ and $\Sp(\mb
H_n)=Q(\mb H_n)$. In particular, $\Sp(\mb H_n)$ is generated by $
s_{u_1}$ and $ s_{v_1}$, where $(u_1,v_1)$ is the standard
hyperbolic pair of $\mb H_n$.
\end{lem}
\bp Note that $(u_1,v_1)$ is an (ordered) $\mb Z_n$-basis for $\mb
H_n$. For any $\varphi\in \Sp(\mb H_n)$, one has
$$\va(u_1)=a_{11}u_1 +a_{21}v_1,\ \va(v_1)=a_{12}u_1+a_{22}v_1$$ where $a_{ij}\in\mb Z_n$. Then with respect to the $\mb Z_n$-basis $(u_1,v_1)$, the matrix of $\va$ is defined to be \be C=\left(
             \begin{array}{cc}
               a_{11} & a_{12} \\
               a_{21} & a_{22}\\
             \end{array}
           \right),\lb{}\ee which is in $\GL(2,\mathbb{Z}_n)$.
 This defines a map $\Sp(\mb H_n)\rt \GL(2,\mb Z_n)$, which is an injective homomorphism. For any $a\in\mb H_n$, $a=i u_1+j v_1$
 with $i,j\in\mb Z_n$. Then the coordinate $\wt a$ of $a$ is denoted $\wt a=(i,j)^t$, the transpose of $(i,j)$. It is clear that \be \wt{\va(a)}=C \wt a.\lb{n}\ee

If we identify $$C_n=\{\omega_n^i|i=0,1,\cdots,n-1\}\rt\mb Z_n,\
\omega_n^i\rt i,$$ then the matrix of the pairing $<,>$ in the $\mb
Z_n$-basis $(u_1,v_1)$ is $J$.

For any $a,b\in\mb H_n$, one has \be <a,b>=\wt a^t J \wt b\lb{m}.\ee

As $\va\in \Sp(\mb H_n)$, by (\ref{n}) and (\ref{m}) one has
$$C^t J C=J,$$ thus $C\in \Sp(2,\mb Z_n)$. So $\Sp(\mb H_n)\st \Sp(2,\mb Z_n)=\SL(2,\mb Z_n)$.

As $$ s_{u_1}(u_1)=u_1,\ s_{u_1}(v_1)=v_1-<u_1,v_1>u_1=v_1-u_1$$ so
the matrix of $ s_{u_1}$ is $$A=\left(
             \begin{array}{cc}
               1 & -1\\
               0 & 1 \\
             \end{array}
           \right) .$$
As $$ s_{v_1}(u_1)=u_1-<v_1,u_1>v_1=u_1+v_1,\ s_{v_1}(v_1)=v_1$$ so
the matrix of $ s_{v_1}$ is $$B=\left(
             \begin{array}{cc}
               1 & 0\\
               1 & 1 \\
             \end{array}
           \right) .$$
It is clear that $A,B$ generate $\SL(2,\mb Z_n)$, thus $\Sp(\mb
H_n)=Q(\mb H_n)= \SL(2,\mb Z_n)$. \ep
\begin{coro}\lb{l}
 For any $a\in \mb H_n$, if $\ord(a)=n$ then $a$ is conjugate to $u_1$ under $\Sp(\mb H_n)$; otherwise $a$ is conjugate to $l u_1$ for
 some $l\in\mb Z_n$.
  In particular all the elements in $\mb H_n$ with order $n$ are
conjugate under $\Sp(\mb H_n)$.
\end{coro}
\bp Assume $a=(i,j)$.  As $\mb Z_n$ is a principal ideal ring, the
ideal $I(i,j)$ in $\mb Z_n$ generated by $i,j$ must be generated by
some $l\in\mb Z_n$. So $I(i,j)=I(l)$.

If $\ord(a)=n$, $I(l)=\mb Z_n$. So there exists $k,m\in \mb Z_n$
such that $ik+jm=1$. Let $A=\left(
             \begin{array}{cc}
               k& m\\
               -j & i \\
             \end{array}
           \right)\in \SL(2,\mb Z_n)$. Then $A(i,j)^t=(1,0)^t$, so $a$ is conjugate to $u_1$ by $A$.

If $\ord(a)<n$, then $(i,j)=l(i^{'},j^{'})$ for some
$a^{'}=(i^{'},j^{'})\in \mb H_n$ as $l$ is the greatest common
devisor of $i,j$. Then $\ord(a^{'})=n$ and there exists some $A\in
\SL(2,\mb Z_n)$ such that $A(a^{'})^t=(1,0)^t$ and $A a^t=(l,0)^t$.
Thus $a$ is conjugate to $l u_1$ by $A$.  \ep

\begin{lem}\lb{o}
 Let $H=\mb H_n\op \mb H_n$ and $\phi:H\rt H,\ (a,b)\mapsto (b,a)$. Then $\phi\in Q(H)$.
\end{lem}
\bp

Let $(u_1,v_1)$ (resp. $(u_2,v_2)$) be the standard hyperbolic pair
in $\mb H_n\op 0$ (resp. $0\op \mb H_n$). Then $\phi$ maps
$u_1$ to $u_2$ and $v_1$ to $v_2$. Let $x=v_1+v_2$, then $(u_1,x)$ and
$(u_2,x)$ are both hyperbolic pairs of order $n$. Let $H_1=Span
(u_1,x)$, then $H=H_1\op H_1^{\perp}$. By Corollary \ref{l}, there
exists $$\tau=(\tau^{'},1)\in \Sp(H_1)\times \Sp(H_1^{\perp})\st
\Sp(H)$$ such that $\tau(x)=u_1$. Similarly there exists $\va\st
\Sp(H)$ such that $\va(u_2)=x$. Then $\tau\va(u_2)=u_1$. Let
$v_2^{'}=\tau\va(v_2)$.

Assume $<v_2^{'},v_1>=\omega_n^i$. Then
$$ s_{u_1}^{i-1}(u_1)=u_1,$$
$$ s_{u_1}^{i-1}(v_2^{'})=v_2^{'}-(i-1)(<u_1,v_2^{'}>\cdot u_1)=v_2^{'}-(i-1)u_1.$$ Let $v_2^{''}=v_2^{'}-(i-1)u_1$. Then $<v_2^{''},v_1>=\omega_n$.
 Let $q=v_2^{''}-v_1$. Note $\ord(q)=n$, then $$ s_q(u_1)=u_1-<q,u_1>\cdot q=u_1,$$
 $$ s_q(v_2^{''})=v_2^{''}-<q,v_2^{''}>\cdot q=v_2^{''}-\omega_n\cdot q=v_2^{''}-q=v_1.$$ So the map $\nu= s_q  s_{u_1}^{i-1}
\tau\va\in Q(H)$ maps $(u_2,v_2)$ to $(u_1,v_1)$. Then $\nu\phi$
fixes $(u_1,v_1)$ and maps its orthocomplement $0\op \mb H_n$
isometrically onto itself. Then there exists
$\theta=(1,\theta^{'})\in \Sp(\mb H_n)\times \Sp(\mb H_n)$ such that
$\theta\nu\phi=1$. So $\phi\in Q(H)$ as $\theta$ and $\nu$ are both
generated by transvections.

\ep
\begin{lem}
 Assume that $H=H(p)$ for some prime $p$. Then by Theorem \ref{d1} (2), $H=\mb H_{p^{r_1}}\op \mb H_{p^{r_2}}\op\cdots \op \mb H_{p^{r_s}}$ for some positive integers
  $r_1,r_2,\cdots,r_s$ with $r_i\geq r_{i+1}$. For any $a=(a_1,\cdots,a_s)\in H$ with order $p^{r_1}$, there exists $\phi\in Q(H)$ such that
   $\phi(a)=b=(b_1,\cdots,b_s)$ with $\ord({b_1})=p^{r_1}$.
\end{lem}
\bp
 In this case one has $\ord(a)=Max_{i=1}^s\{\ord(a_i)\}$. If $\ord(a_1)=p^{r_1}$ then we take $\phi=1$. If $\ord(a_i)=p^{r_1}$ for some $i>1$, then $r_i=r_1$.
 Then as the subgroups $\mb H_{p^{r_1}}$ and $\mb H_{p^{r_i}}$ of $H$ are isometric, by Lemma \ref{o} there exists some
 $\phi\in Q(\mb H_{p^{r_1}} \op \mb H_{p^{r_i}})\st Q(H)$
  such that $$\phi(a)=\phi(a_1,\cdots,a_i,\cdots,a_s)=(a_i,\cdots,a_1,\cdots,a_s)=b. $$ Then $b_1=a_i$ and $\ord({b_1})=p^{r_1}$.
\ep

\begin{lem}
 Assume $H=\mb H_{l_1}\op \mb H_{l_2}\op \cdots \op \mb H_{l_k}$ with $l_i|l_{i-1}$ for all $i$. Then for any $a=(a_1,\cdots,a_k)\in H$ with order $l_1$,
 there exists $\phi\in Q(H)$ such that $\phi(a)=b=(b_1,\cdots,b_k)$ with $\ord(b_1)=l_1$.
\end{lem}
\bp Let $p_1,\cdots,p_s$ be the set of primes dividing $|H|$. Then
$H=\op_i H(p_i)$ and $H(p_i)=\mb H_{l_1}(p_i)\op \mb H_{l_2}(p_i)\op
\cdots \op \mb H_{l_k}(p_i)$. Let $\pi_i:H\rt H(p_i)$ be the
projection. Then $$\pi_i(a)=(a_{1i},a_{2i},\cdots,a_{ki})\in \mb
H_{l_1}(p_i)\op \mb H_{l_2}(p_i)\op \cdots \op \mb H_{l_k}(p_i).$$
One has $a=\sum_i \pi_i(a)$. By last lemma there exists $\phi_i\in
Q(H(p_i))\st Q(H)$ such that
$\phi_i(\pi_i(a))=(b_{1i},b_{2i},\cdots,b_{ki})$ with
$\ord(b_{1i})=\ord(\pi_i(a))$. Let $\phi=\Pi_i \phi_i$. Then
$\phi(a)=b=(b_1,b_2,\cdots,b_k)\in \mb H_{l_1}\op \mb H_{l_2}\op
\cdots \op \mb H_{l_k}$, where $b_1=\sum_i b_{1i}$ and
$\ord(b_{1})=\ord(b)=\ord(a)=l_1$.

\ep

\begin{lem}
 Assume $H=\mb H_{l_1}\op \mb H_{l_2}\op \cdots \op \mb H_{l_k}$ with $l_i|l_{i-1}$ for all $i$. Then for any $a=(a_1,a_2,\cdots,a_k)\in H$ with order $l_1$,
 there exists $\phi\in Q(H)$ such that $\phi(a)=b=(u_1,0,\cdots,0)$, where $u_1=(1,0)\in \mb H_{l_1}$.
\end{lem}

\bp By last lemma we can assume that $\ord({a_1})=l_1$. Then use
induction on the number $t$ of nonzero elements in
$\{a_1,\cdots,a_k\}$. The case $t=1$ follows from Corollary \ref{l}.

Assume there are $t=l\geq 2$ nonzero elements in
$\{a_1,a_2,\cdots,a_k\}$ and the result holds for $l-1$. Without
loss of generality we can assume $a_2\neq 0$. As $\ord({a_1})=l_1$,
there exists some $b_1\in \mb H_{l_1}\st H$ such that
$<b_1,a_1>=\omega_{l_1}$. Let $b=(b_1,a_2,0\cdots,0)$. Then
$\ord({b})=l_1$ and $<b,a>=\omega_{l_1}$, so \be \begin{split}
 s_b(a)&=a-<b,a>\cdot b\\&=a-\omega_{l_1}\cdot b=a-b
\\&=(a_1-b_1,0,a_3,\cdots,a_k).
\end{split}\ee So $ s_b(a)=(a_1-b_1,0,a_3,\cdots,a_k)$ differs with $a$
only in the first and second position. As $ s_b(a)$ has order $l_1$
and has $l-1$ nonzero elements, by induction  there exists
$\phi_1\in Q(H)$ such that $\phi_1( s_b(a))=(u_1,0,\cdots,0)$. Then
$\phi=\phi_1  s_b\in Q(H)$ has the desired property and the result
holds for $t=l$.

\ep
\begin{coro}\lb{j}
 $Q(H)$ acts transitively on the set of elements in $H$ with maximal order.
\end{coro}

\begin{lem}
  Assume $H=\mb H_{l_1}\op \mb H_{l_2}\op \cdots \op \mb H_{l_k}$ with $l_i|l_{i-1}$ for all $i$ and $G=\Sp(H)$. For $i=1,\cdots,k$ let $(u_i,v_i)$ be the standard hyperbolic pair in $\mb H_{l_i}$. Then (1) $G_{u_1}=Q(H)_{u_1}$. (2) $G=Q(H)$.
\end{lem}
\bp We will prove them by induction on $k$.

The case $k=1$ follows from Lemma \ref{h} as we proved there
$\Sp(H)=Q(H)$ if $H=\mb H_{l_1}$. Let $k>1$. Assume (1) and (2) holds for
$k-1$.

Assume $\si\in G_{u_1}$. As
$<u_1,v_1>=<\si(u_1),\si(v_1)>=<u_1,\si(v_1)>$,
$<u_1,\si(v_1)-v_1>=0$ so
$$\si(v_1)=v_1+j_1 u_1+\sum_{i=2}^k (p_i u_i+ q_i v_i)$$ for some
$j_1,p_i,q_i\in \mb Z$. By Corollary \ref{l} there exists $\phi_i\in
Q(\mb H_{l_i})\st Q(H)$ such that $\phi_i(p_i u_i+ q_i v_i)=j_i u_i$
for $i\geq 2$. Let $\phi=\Pi_{i=2}^k\ \phi_i$. Then  $\phi\in
Q(H)_{u_1}$ and $$\phi \si(v_1)=v_1+j_1 u_1+\sum_{i=2}^k j_i u_i.$$

For any $i$ with $2\le i\le k$, $ s_{u_1+u_i}(u_t)=u_t$ for
$t=1,\cdots,k$. As $\ord(u_1+u_i)=l_1$, \be
\begin{split} s_{u_1+u_i}(v_1)&=v_1-<u_1+u_i,v_1>\cdot (u_1+u_i)\\&=v_1-\omega_{l_1}\cdot (u_1+u_i)\\&=v_1-(u_1+u_i).\end{split}\ee

Let $\tau=\Pi_{i=2}^k\  s_{u_1+u_i}^{j_i}$. Then $\tau\in
Q(H)_{u_1}$ and
$$ \tau(\phi\si(v_1)
)=v_1+(j_1-\sum_{i=2}^k j_i)u_1.$$
So $\tau\phi\si$ preserves $\mb H_{l_1}$ and also $\mb
H_{l_1}^{\perp}=\mb H_{l_2}\op \cdots \op \mb H_{l_k}$ thus
$$\tau\phi\si\in \Sp(\mb H_{l_1})\times \Sp(\mb H_{l_2}\op \cdots
\op \mb H_{l_k}).$$ By induction $$\Sp(\mb H_{l_1})\times \Sp(\mb
H_{l_2}\op \cdots \op \mb H_{l_k})=Q(\mb H_{l_1})\times Q(\mb
H_{l_2}\op \cdots \op \mb H_{l_k}),$$ so $$\tau\phi\si\in Q(\mb
H_{l_1})\times Q(\mb H_{l_2}\op \cdots \op \mb H_{l_k})\st Q(H).$$
As $\tau,\phi\in Q(H)$, one also has $\si\in Q(H)$. So $\si\in
Q(H)_{u_1}$ then $G_{u_1}\st Q(H)_{u_1}$. As $G\supset Q(H)$, one
must have $G_{u_1}= Q(H)_{u_1}$. So (1) holds for $k$.

By Corollary \ref{j}, $Q(H)$ acts transitively on the set of
elements in $H$ with maximal order, so does $G$. As
$$G_{u_1}=Q(H)_{u_1}\ \an \ |G/G_{u_1}|=|Q(H)/Q(H)_{u_1}|,$$ so
$|G|=|Q(H)|$ thus $G=Q(H)$. So (2) also holds for $k$. \ep

Now we have proved the following theorem.
\begin{theorem}\lb{c3}
Let $H$ be a finite nonsingular symplectic abelian group. Then
$\Sp(H)$ is generated by the set of transvections on $H$.
\end{theorem}
\begin{rem}
If $H=\mb H_n^k$, the $k$-fold direct sum of $\mb H_n$, then by
choosing some suitable $\mb Z_n$-basis of $\mb H_n$, it is easy to
see that
$$\Sp(\mb H_n^k)\cong\Sp(2k,\mb Z_n),$$ where
$$\Sp(2k,\mathbb{Z}_n)=\{A\in GL(2k,\mathbb{Z}_n)|A^t J_{2k}
A=J_{2k}\}$$ with $$J_{2k}=\left(
             \begin{array}{ccccc}
               0 & 1 & & &\\
               -1 & 0 & & &\\
                  &  & \ddots & &\\
                & & & 0 & 1 \\
               & & &-1 & 0 \\
             \end{array}
           \right)\in M(2k,\mathbb{Z}_n).$$ Thus this theorem implies
           that in particular $\Sp(2k,\mb Z_n)$ is generated by the
           transvections.
\end{rem}
\section{Fine gradings of Lie algebras
and their Weyl groups }
 \setcounter{equation}{0}\setcounter{theorem}{0}
4.1. The details of this subsection can be found in Section 4 of \cite{hg}. We include it for completeness.

We always assume that $L$ is a complex simple Lie algebra. Let
$\Aut(L)$ be its automorphism group and $\Int(L)$ the identity
component of $\Aut(L)$, called the inner automorphism group
 of $L$. It is clear that $\Aut(L)$ and $\Int(L)$ are both algebraic groups.

 Let $\Lambda$ be a additive  abelian group. A
$\Lambda$-\textit{grading} $\Gamma$ on $L$ is the decomposition of
$L$ into direct sum of subspaces $$\Gamma: L=\oplus_{\ga\in \Lambda}
\ L_\ga$$ such that
$$[L_\ga,L_\dt]\st L_{\ga+\dt},\ \forall\\ \ga,\dt\in \Lambda.$$ Let
$\ddt=\{\ga\in \Lambda|L_\ga\neq 0\}$. We will always
assume that $\Lambda$ is generated by $\ddt$, otherwise it could be
replaced by its subgroup generated by $\ddt$. So $\Lambda$ is always finitely generated.

Given a $\Lambda$-grading $\Gamma$ on $L$, let $$K=\widehat
\Lambda=_{def} \Hom(\Lambda, \mathbb{C}^{\times})$$ be the abelian
group of characters of $\Lambda$. Then $K$ acts on $L$ by
$$\si\cdot X=\si(\ga) X,\ \forall\\ X\in L_\ga,\ \forall\ \ga\in \Lambda,\
\forall\\si\in K.$$ This defines an injective homomorphism $K\rt
\Aut(L)$. So $K$ can be viewed as a subgroup of $\Aut(L)$. Recall
that an algebraic group is called \textit{diagonalizable} if it is
abelian and consists of semisimple elements. It is easy to see that
$K$ is a diagonalizable algebraic subgroup of $\Aut(L)$.

Conversely, given a diagonalizable algebraic subgroup $K$ of
$\Aut(L)$, let $$\Lambda=\widehat K=_{def} \Hom(K,
\mathbb{C}^{\times})$$ be the (additive) abelian group of
homomorphisms from $K$ to $\mb C^{\times}$ as algebraic groups. Then
one has a $\Lambda$-grading on $L$: $$\Ga:\ \ L=\oplus_{\gamma\in
\Lambda}\ L_\ga,$$ where $L_\ga=\{X\in L|\si \cdot X=\ga(\si) X,\
\forall\ g\in K\}.$ Let $$\ddt=\ddt(L,K)=_{def} \{\ga\in
\Lambda|L_{\ga}\neq 0\}.$$ We call $\ddt$ the set of \textit{roots}
of $K$ in $L$.

Thus there is a natural one-to-one correspondence between gradings of
$L$ by finitely generated abelian groups and diagonalizable
algebraic subgroups of $\Aut(L)$.  A grading is called \textit{inner} if the respective diagonalizable subgroup is in
 $\Int(\L)$. A grading (resp. inner grading) of $L$ is called fine if it could not be further refined by any other grading (resp. inner grading).
It is clear that the bigger the
diagonalizable algebraic subgroup $K$ is, the finer the
corresponding grading is. Thus a grading (resp. inner grading) $\Ga$ on $L$ is
fine if and only if the corresponding diagonalizable subgroup $K$
is a maximal diagonalizable subgroup of $\Aut(L)$ (resp. $\Int(L)$).

 Let $G$ be either $\Aut(L)$ or $\Int(L)$. Let $K$ be a {maximal} diagonalizable subgroup of $G$,
$\Ga$ the grading on $L$ induced by the action of $K$. One could
define the Weyl group $W_G(\Ga)$ of the grading $\Ga$ with respect
to $G$, see Definition 2.3 of \cite{hg}, to describe the symmetry of
the grading $\Ga$.

One has the following result in \cite{hg}.
\begin{prop}[Corollary 2.6 of \cite{hg}]
 Let $L$ be a simple Lie algebra and $G=\Int(L)$. Let $K$ be a maximal diagonalizable subgroup of $G$ and $\Ga$ be the corresponding grading
 on $L$ induced by the action of $K$. Let $W_G(K)=N_G(K)/K$ be the Weyl group
 of $K$ with respect to $G$. Then one has $W_G(\Ga)=W_G(K)$. \end{prop}
 \bigskip
4.2. From now on we will always assume $K\st G=\Int(L)$ to be a {finite
maximal diagonalizable subgroup} and $\ddt=\ddt(L,K)$. Let $B$ be the Killing form on
$L$. Recall that a linear subspace $S$ of $L$ is called a toral
subalgebra if $[S,S]=0$ and the endomorphism $\ad_X$ is semisimple
for each $X\in S$.

As $L$ is simple, the adjoint map
$$\ad:L\rt \ad(L),\ X\mapsto \ad_X$$ is a $G$-equivariant isomorphism.
\begin{defi} For any
$\ga\in \ddt$, let
$$L_{[\ga]}=_{def}\oplus_{k\in \mathbb{Z}}\ L_{k\ga}.$$ \end{defi}

\begin{prop}\lb{d0}
(1) One has $L_0=0$, i.e., $0\notin \ddt$.

(2) Assume $\ga,\dt\in\ddt$ and $\ga+\dt\neq 0$. Then $B|L_\ga\times
L_\dt=0$. For any $X\in L_\ga$ with $X\neq 0$, there exists $Y\in
L_{-\ga}$ such that $B(X,Y)\neq 0$.

(3) For any $\ga\in \ddt$, $L^{\Ker \ga}=L_{[\ga]}$ and is a toral
subalgebra of $L$. One has $Lie\ Z(Ker \ga)_0=ad(L_{[\ga]})$ and
$Z(Ker \ga)_0$ is an algebraic torus (isomorphic to some
$(\mathbb{C}^\times)^i$).

\end{prop}
\bp (1) As $K$ is a maximal diagonalizable subgroup, $Z_G(K)=K$ by Lemma 2.2 of \cite{hg}. As
$K$ is finite, $$\ad(L^K)=\ad(L)^K=\Lie Z_G(K)=\Lie K=0.$$  So
$L_0=L^K=0$.

(2) Assume $\ga+\dt\neq 0$. For any $X\in L_\ga$ and $Y\in L_\dt$,
$\ad_X \ad_Y$ maps each $L_\zeta$ into $L_{\zeta+\ga+\dt}$ thus
$B(X,Y)=0$. Then $B|L_\ga\times L_\dt=0$. Because $B$ is nonsingular
on $L$, the second statement then follows from the first one.

(3) We first prove $L^{\Ker \ga}=L_{[\ga]}$. Choose some $\si\in K$
satisfying $\ga(\si)=\omega_m$, where $m$ is the order of $\ga$.
Then $\ga(\si)$ is a generator of the cyclic group $\ga(K)\cong
K/\Ker \ga$. $L^{\Ker \ga}$ is the direct sum of those
 $L_\bt$ with $\bt$
being identity on $\Ker \ga$. Then $\bt(\si)=\omega_m^k$ for some
integer $k$ as $\si^m\in \Ker \ga $. Then
$\bt(\si)=\ga(\si)^{k}=(k\ga)(\si)$. As $\Ker \ga$ and $\si$
generate $K$, $\bt=k\ga$. Thus $L^{\Ker \ga}=\oplus_{k}
L_{k\ga}=L_{[\ga]}$.

One has $[L_0,L_{[\ga]}]=0$ as $L_0=0$. Then by Proposition 3.6 of
\cite{h}, $L_{[\ga]}$ is a toral subalgebra of $L$. As
$\Lie\Int(L)=\ad(L)$, it is clear that
$$\Lie\ Z(\Ker \ga)_0=\Lie\ Z(\Ker \ga)=\ad(L)^{\Ker \ga}=\ad(L^{\Ker
\ga})=\ad(L_{[\ga]}).$$Thus $Z(Ker \ga)_0$ is an algebraic torus.

\ep

\begin{rem}If $L$ is a semisimple Lie algebra then all the
results in this section still hold.
\end{rem}

\section{Finite maximal diagonalizable subgroups of $\PGL(n,\mb C)$ and anti-symmetric pairings on them}
 \setcounter{equation}{0}\setcounter{theorem}{0}

 Let $n\in \mathbb{Z}_+$. Recall $\omega_n=e^{2\pi i/n}$. Let
$$Q_n=diag(1,\omega_n,\omega_n^2,\cdots,\omega_n^{n-1})$$ and  $$P_n=\left(
             \begin{array}{ccccc}
               0 & 1 & 0 &\cdots &  0 \\
               0 & 0 & 1 &\cdots &  0 \\
               \vdots& \vdots& \vdots & \ddots  & \vdots  \\
               0 & 0 & 0 &\cdots &  1 \\
               1 & 0 & 0 &\cdots &  0\\
             \end{array}
           \right).$$ Let $\Pi_n=\{\omega_n^j P_n^k Q_n^l\ |j,k,l=0,1,\cdots,n-1\}$. Note $P_n Q_n=\omega_n Q_n P_n$. This is a subgroup of $\GL(n,\mathbb{C})$,
           called the \textit{Pauli group} of rank $n$.

Let $D_n$ be the subgroup of diagonal matrices
of $\GL(n,\mb{C})$. Let $\texttt{P}_n$ and $\texttt{D}_n$
 be the respective images of $\Pi_n$ and $D_n$ under the adjoint action on $M(n,\mb C)$. One knows that
$$\texttt{P}_n=\{Ad_{P_n}^iAd_{Q_n}^j|i,j=0,\cdots,n-1\}\cong \mb
Z_n\times \mb Z_n,$$ and that
 $\texttt{D}_n$ and $\texttt{P}_n$ are both maximal
diagonalizable subgroups of $\PGL(n,\mathbb{C})$.

Let $L=sl(n,\mathbb{C})$ and $G=\Int(L)\cong \PGL(n,\mathbb{C})$.
One has the standard isomorphism
$$M=M(t,\mathbb{C})\ot M(l_1,\mathbb{C})\ot \cdots \ot M(l_k,\mathbb{C})\rt M(n,\mathbb{C}),$$ where $n=t l_1\cdots l_k$ and $l_i|l_{i-1}$ for all $i$. It induces injective homomorphisms \be
D_t\ot\Pi_{l_1}\ot\cdots\ot \Pi_{l_k}\st S=\GL(t,\mathbb{C})\ot
\GL(l_1,\mb C)\ot \cdots \ot \GL(l_k,\mathbb{C})\rt
\GL(n,\mathbb{C}).\lb{a8}\ee

   Let $A=A_0\ot A_1\ot\cdots \ot A_k\in S$. Then for $X=X_0\ot X_1\ot\cdots\ot X_k\in M$,
 $$\Ad_A(X)=\Ad_{A_0}(X_0)\ot\Ad_{A_1}(X_1)\ot\cdots\ot \Ad_{A_k}(X_k).$$
 Thus the adjoint action induces homomorphisms $$\GL(n,\mathbb{C})\rt \PGL(n,\mathbb{C}),$$
 \bee \GL(t,\mathbb{C})\ot
\GL(l_1,\mb C)\ot \cdots \ot \GL(l_k,\mathbb{C})\rt \PGL(t,\mathbb{C})\times \PGL(l_1,\mb C)\times \cdots \times
 \PGL(l_k,\mathbb{C}),\eee
 \be\lb{a4} A=A_0\ot A_1\ot\cdots \ot A_k\mapsto Ad_A=(\Ad_{A_0},\Ad_{A_1},\cdots,\Ad_{A_k})\ee and by restriction
 $$D_t\ot\Pi_{l_1}\ot\cdots\ot \Pi_{l_k}\rt  \texttt{D}_t\times \texttt{P}_{l1}\times\cdots\times\texttt{P}_{l_k}.$$
Thus by (\ref{a8}) one has injective homomorphisms

\be \CD \texttt{D}_t\times
\texttt{P}_{l1}\times\cdots\times\texttt{P}_{l_k}\st
\PGL(t,\mathbb{C})\times \PGL(l_1,\mb C)\times \cdots \times
\PGL(l_k,\mathbb{C}) @>\phi>> \PGL(n,\mathbb{C}).
\endCD\lb{a9}
\ee

We will identify $\texttt{D}_t\times \texttt{P}_{l_1}\times \cdots
\times \texttt{P}_{l_k}$ with its image in $\PGL(n,\mb{C}).$
 \begin{theorem}[Theorem 3.2 of \cite{hpp}] Any maximal diagonalizable subgroup of $\PGL(n,\mb{C})$ is conjugate to one and only one of the $\texttt{D}_t
 \times \texttt{P}_{l_1}\times \cdots \times \texttt{P}_{l_k}$ with $n=t l_1\cdots l_k$ and each $l_i$ dividing $l_{i-1}$.
\end{theorem}

 \begin{coro} \lb{b2} Any finite maximal diagonalizable subgroup of $\PGL(n,\mb{C})$ is conjugate to one and only one of the  $\texttt{P}_{l_1}\times \cdots \times \texttt{P}_{l_k}$
 with $n=l_1\cdots l_k$ and each $l_i$ dividing $l_{i-1}$.
\end{coro}

In the case $K=\texttt{P}_{l_1}\times \cdots \times
\texttt{P}_{l_k}$ is a finite maximal diagonalizable subgroup of
$\PGL(n,\mb{C})$, (\ref{a9}) becomes
\be \CD
\texttt{P}_{l1}\times\cdots\times\texttt{P}_{l_k}\st
\PGL(l_1,\mb C)\times \cdots \times
\PGL(l_k,\mathbb{C}) @>\phi>> \PGL(n,\mathbb{C})
\endCD\lb{b0}
\ee
\bigskip

Let $K\st \PGL(n,\mb{C})$ be a finite maximal diagonalizable
subgroup. Let
$$p:\GL(n,\mathbb{C})\rt \PGL(n,\mathbb{C})$$ be the projection.

\begin{defi} \lb{d} For
any $\si\in K$ fix some $\widetilde{\si}\in p^{-1}(\si)$. For any
$\si,\tau\in K$, $\widetilde{\si} \widetilde{\tau}
\widetilde{\si}^{-1} \widetilde{\tau}^{-1}=lI_n$ as
$p(\widetilde{\si} \widetilde{\tau} \widetilde{\si}^{-1}
\widetilde{\tau}^{-1})=1$. Clearly $l$ is independent of the
preimages  $\widetilde{\si},\widetilde{\tau}$ chosen. Define
$<\si,\tau>= l$.
\end{defi}
\begin{lem}[Lemma 3.4 of \cite{hg}] The map $<,>:K\times K\rt \mathbb{C}^{\times}$ is an anti-symmetric pairing on $K$, which is invariant under $N_G(K)$.
\end{lem}

\begin{prop}\lb{b1}
Let $G=\PGL(n,\mb{C})$ and $K$ be a maximal diagonalizable subgroup
of $G$. If $K=\texttt{P}_{l_1}\times \cdots \times \texttt{P}_{l_k}$
with $n=l_1\cdots l_k$ and each $l_i$ dividing $l_{i-1}$, then $<,>$
is nonsingular on $K$. Thus $(K,<,>)$ is a nonsingular symplectic
abelian group isometric to $\mb H_{l_1}\op \cdots \op \mb H_{l_k}$.
\end{prop}
\bp For $i=1,\cdots,k$, let
$$\si_i=(1,\cdots,1,Ad_{P_{l_i}},1,\cdots,1)\in K,
\ \tau_i=(1,\cdots,1,Ad_{Q_{l_i}},1,\cdots,1)\in K$$ where
$Ad_{P_{l_i}}$ and $Ad_{Q_{l_i}}$ are in the $i$-th position. Then
$\{\si_i,\tau_i|i=1,\cdots,k\}$ is a set of generators of $K$ and
any element in $K$ can be written uniquely as
$\si_1^{i_1}\tau_1^{j_1}\cdots \si_k^{i_k}\tau_k^{j_k}$. By simple
computation one has for $i \neq j$,
$$<\si_i,\tau_i>=\omega_{l_i},\ <\si_i,\tau_j>=1,$$  $$<\si_i,\si_j>=1,\ <\tau_i,\tau_j>=1.$$
Thus $(\si_i,\tau_i)$ is a hyperbolic pair of order $l_i$ and spans
the symplectic subgroup $\texttt{P}_{l_i}$ isometric to $\mb
H_{l_i}$. It is clear that such subgroups $\texttt{P}_{l_i}$ are
mutually orthogonal to each other. The map
\be\lb{c4}\texttt{P}_{l_1}\times \cdots \times \texttt{P}_{l_k}\rt
\mb H_{l_1}\op \cdots \op \mb H_{l_k},\
\si_1^{i_1}\tau_1^{j_1}\cdots \si_k^{i_k}\tau_k^{j_k}\mapsto
((i_1,j_1),\cdots,(i_k,j_k))\ee is clearly an isometry of
nonsingular symplectic abelian groups. \ep

As a corollary of Corollary \ref{b2}, Proposition \ref{b1} and Corollary
\ref{d9} one has the following result.
\begin{prop}\lb{d6}
There is a one-to-one correspondence between conjugacy class of
finite maximal diagonalizable subgroups of $\PGL(n,\mathbb{C})$ and
nonsingular symplectic abelian groups of order $n^2$ .
\end{prop}

\section{Weyl groups of finite maximal diagonalizable subgroups of $\PGL(n,\mb C)$ }
 \setcounter{equation}{0}\setcounter{theorem}{0}
Recall that $L=sl(n,\mb C)$ and $K$ is a finite maximal
diagonalizable subgroup of $G=\PGL(n,\mb{C})$. First we will
describe the grading of $sl(n,\mb C)$ and $gl(n,\mb C)$ induced by
the action of $K$.

At first let $K=\texttt{P}_n$. The character group $\wh{K}$, written
additively, is generated by $\bt_n$ and $\al_n$, which are dual to
$\Ad_{P_n},\Ad_{Q_n}$:
$$\bt_n(\Ad_{P_n})=\omega_n,\ \bt_n(\Ad_{Q_n})=1,$$
$$\al_n(\Ad_{P_n})=1,\ \al_n(\Ad_{Q_n})=\omega_n.$$
Thus $\wh{K}=\{i\bt_n+j\al_n|(i,j)\in \mathbb{Z}_n\times
\mathbb{Z}_n\}\cong \mathbb{Z}_n^2$.

One has
$$\bt_n(\Ad_{P_n}^{i}\Ad_{Q_n}^{j})=\omega_n^i=<\Ad_{Q_n}^{-1},\Ad_{P_n}^{i}\Ad_{Q_n}^{j}>$$
and
$$ \al_n(\Ad_{P_n}^{i}\Ad_{Q_n}^{j})=\omega_n^j=<\Ad_{P_n},\Ad_{P_n}^{i}\Ad_{Q_n}^{j}>,$$
 so $$\bt_n^*=\Ad_{Q_n}^{-1},\ \al_n^*=\Ad_{P_n}.$$ As by (\ref{b3}) one has $(\ga+\dt)^*=\ga^*\dt^*$,
\be \lb{a1}(i\bt_n+j\al_n)^*=\Ad_{P_n}^{j}\Ad_{Q_n}^{-i}.\ee As
$$\Ad_{Q_n}(Q_n^i P_n^j)=\omega_n^{-j}Q_n^i P_n^j=(i\bt_n-j\al_n)(\Ad_{Q_n})Q_n^i P_n^j,$$
$$\Ad_{P_n}(Q_n^i P_n^j)=\omega_n^{i}Q_n^i P_n^j=(i\bt_n-j\al_n)(\Ad_{P_n})Q_n^i P_n^j,$$ one has
\be \lb{z} Q_n^i P_n^j\in L_{i\bt_n-j\al_n}.\ee In particular
$$P_n\in L_{-\al_n},\ Q_n\in L_{\bt_n}.$$ Note that $tr (Q_n^i P_n^j)=0$ unless
$(i,j)=(0,0)$. Let \be X_{i\bt_n-j\al_n}=Q_n^i P_n^{j}.\lb{a2}\ee
Then one has the following gradings

$$gl(n,\mathbb{C})=\oplus_{(i,j)}\ \mathbb{C} Q_n^i P_n^j=\op_{\ga\in \wh K}\ \mb C X_{\ga},$$
$$sl(n,\mathbb{C})=\oplus_{(i,j)\neq (0,0)}\ \mathbb{C} Q_n^i P_n^j=\op_{\ga\neq 0}\ \mb C X_{\ga}.$$ So
$\ddt(gl(n,\mathbb{C}),K)=\wh{K}$ and
$\ddt(sl(n,\mathbb{C}),K)=\wh{K}\setminus \{0\}$. Note that each
root space is one-dimensional, and for any $\ga\in\wh K$, by
(\ref{a1}) and (\ref{a2}) one has \be
\ga^*=(Ad_{X_\ga})^{-1}.\lb{a3}\ee The following result is
originally Theorem 10 of \cite{jpt}, for its proof we refer the
readers to Proposition 4.4 of \cite{hg}.
\begin{theorem}
Let $G=\PGL(n,\mb C)$ and $K=\texttt{P}_n$. One has $W_G(K)\cong
\SL(2,\mb Z_n)$ and is generated by $s_{\al_n}$ and $s_{\bt_n}$.
\end{theorem}

Next let $K=\texttt{P}_{l_1}\times \cdots \times \texttt{P}_{l_k}$
with $n=l_1\cdots l_k$ and each $l_i$ dividing $l_{i-1}$. As
$M(n,\mb C)=M(l_1,\mb C)\ot\cdots\ot M(l_k,\mb C)$, and $M(l_i,\mb
C)=\op_{\ga\in \wh{\texttt{P}_{l_1}}}\mb C X_\ga$, one has
$$ M(n,\mb C)=\op_{(\ga_1,\cdots,\ga_k)}\ \mb C X_{\ga_1}\ot\cdots\ot X_{\ga_k}.$$

Note that $\wh K=\wh{\texttt{P}_{l_1}}\times \cdots \times
\wh{\texttt{P}_{l_k}}$. Let $\ga=(\ga_1,\cdots,\ga_k)\in \wh K$ with
$\ga_i\in \wh{\texttt{P}_{l_i}}$. For any
$\si=(\si_1,\cdots,\si_k)\in K$, $\ga(\si)=\ga_1(\si_1)\cdots
\ga_k(\si_k)$ and one has \bee
\begin{split} \si\cdot X_{\ga_1}\ot\cdots\ot
X_{\ga_k}&=\si_1\cdot X_{\ga_1}\ot\cdots\ot \si_k\cdot
X_{\ga_k}\\&=\ga_1(\si_1)X_{\ga_1}\ot\cdots\ot
\ga_k(\si_k)X_{\ga_k}\\&=\ga(\si) X_{\ga_1}\ot\cdots\ot
X_{\ga_k}.\end{split}\eee

 So $X_{\ga_1}\ot\cdots\ot
X_{\ga_k}\in L_\ga$. Note that $tr (X_{\ga_1}\ot\cdots\ot
X_{\ga_k})=\Pi_i tr(X_{\ga_i})$, which is nonzero if and only if
$\ga_1=\cdots=\ga_k=0$. Let
$$Y_\ga=X_{\ga_1}\ot\cdots\ot X_{\ga_k}.$$ Then
$$gl(n,\mathbb{C})=\op_{\ga\in \wh K}\ \mb C Y_\ga$$ and
$$sl(n,\mathbb{C})=\op_{\ga\neq 0}\ \mb C Y_\ga.$$

So $\ddt(gl(n,\mb C),K)=\wh K$ and $\ddt(sl(n,\mb C),K)=\wh
K\setminus \{0\}$. Note that each root space is also one-dimensional
and consists of semisimple elements.

\begin{lem}For any $\ga\in \wh K$, $\ga^*=(\Ad_{Y_\ga})^{-1}$. \end{lem}

\bp For any $Ad_X\in K$, as $Y_\ga$ is invertible, $$Y_\ga^{-1} X
Y_\ga X^{-1}=Y_\ga^{-1} (\ga(Ad_X)Y_\ga)=\ga(Ad_X) I,$$ so
$<(\Ad_{Y_\ga})^{-1},Ad_X>=\ga(Ad_X)$. Thus
$\ga^*=(\Ad_{Y_\ga})^{-1}$.
 \ep

  Recall in (\ref{b0}) one has the
embedding \bee \CD
\texttt{P}_{l1}\times\cdots\times\texttt{P}_{l_k}\st \PGL(l_1,\mb
C)\times \cdots \times \PGL(l_k,\mathbb{C}) @>\phi>>
\PGL(n,\mathbb{C}).
\endCD
\eee
 Let $N(\texttt{P}_{l_i})$ be
the normalizer of $\texttt{P}_{l_i}$ in $\PGL(l_i,\mathbb{C})$, then
clearly $\phi$ restricts to \bee \CD
N(\texttt{P}_{l_1})\times\cdots\times N(\texttt{P}_{l_k}) @>\phi>>
\PGL(n,\mathbb{C}).
\endCD
\eee The left hand side is in $N_G(K)$. As
$$N(\texttt{P}_{l_i})/\texttt{P}_{l_i}\cong \SL(2,\mathbb{Z}_{l_i}),$$ one has
\be\lb{s} \SL(2,\mathbb{Z}_{l_1})\times\cdots\times
\SL(2,\mathbb{Z}_{l_k})\st W_G(K).\ee

 Let $\ga\in \ddt(sl(n,\mb C),K)$ and $G=\PGL(n,\mb C)$. Assume the order of $\ga$ is $m$ and choose $\si\in K$ satisfying $\ga(\si)=\omega_m$. As $\si\in Z(Ker \ga)$,
$\Ad_\si$ maps $Z(Ker \ga)_0$ into $Z(Ker \ga)_0$. Let $f_\si:G\rt
G,\ \eta\mapsto\si\eta\si^{-1}\eta^{-1}$. Then $f_\si(Z(Ker
\ga)_0)\st Z(Ker \ga)_0$. Denote
 $Z(\Ker \ga)_0$ by $Z_0$.
\begin{lem}\lb{a6}
 (1) The map $f_\si:Z_0\rt  Z_0,\ \eta\mapsto\si\eta\si^{-1}\eta^{-1}$ is a
continuous epimorphism.

(2) Assume $\ga^*\in Z_0$. Then there exists $\zeta\in Z_0$ with
$f_\si(\zeta)=\ga^*$. One has $\zeta\in N_G(K)$ and
 $Ad_\zeta:K\rt K,\ \tau\mapsto \zeta\tau\zeta^{-1}$ is just the transvection $$ s_\ga:K\rt K,\ \tau\mapsto \tau(\ga^*)^{-\ga(\tau)}$$ as in
 (\ref{b9}). (Note that as a subgroup of G,
 $K$ is a multiplicative abelian group.) Thus $ s_\ga\in W_G(K)$.
\end{lem}
\bp As it was shown in Proposition \ref{d0} (3) that $L_{[\ga]}$ is
a toral subalgebra of $L$, the lemma follows from Lemma 3.7 of
\cite{hg}. \ep

\begin{lem}\lb{c1}  For each $\ga\in\ddt(sl(n,C),K)$, $ s_\ga\in W_G(K)$. \end{lem}
\bp Assume
$\ga=(a_1\bt_{l_1}+b_1\al_{l_1},\cdots,a_k\bt_{l_k}+b_k\al_{l_k})$,
then $$Y_\ga =Q_{l_1}^{a_1}P_{l_1}^{-b_1}\ot\cdots\ot
Q_{l_k}^{a_k}P_{l_k}^{-b_k}.$$

By Corollary \ref{l} and (\ref{s}) there exists some
$Y_\dt=Q_{l_1}^{c_1}\ot\cdots\ot Q_{l_k}^{c_k}$ such that
$Ad_{Y_\ga}$ is conjugate to $Ad_{Y_\dt}$ under $N_G(K)$. Assume $
Y_\dt=Q_{l_1}^{c_1}\ot\cdots\ot Q_{l_k}^{c_k}$ as an element of
$\GL(n,\mathbb{C})$ has order $m$. Then $L_{i\dt}=\mb C Y_{\dt}^i$
for $i=1,\cdots,m-1$ and $L_{[\dt]}=\oplus_{i=1}^{m-1} \mb C
Y_\dt^i$ is an abelian Lie algebra consisting of semisimple
elements. We will show $\Ad_{Y_\dt}\in Z(\Ker \dt)_0$, then
$\Ad_{Y_\ga}\in Z(\Ker \ga)_0$ as $Ad_{Y_\ga}$ and $Ad_{Y_\dt}$ are
conjugate. Thus $\ga^*=(\Ad_{Y_\ga})^{-1}\in Z(\Ker \ga)_0$ and $
s_\ga\in W_G(K)$ by Lemma \ref{a6}.

 The set $D_i$ of
eigenvalues of $Q_{l_i}^{c_i}$ is a cyclic group for each $i$. Let
$D$ be the set of eigenvalues of $Y_\dt$. For any $a,b\in D$,
$a=a_1\cdots a_k,\ b=b_1\cdots b_k$ with $a_i,b_i\in D_i$. Then
$ab^{-1}=(a_1 b_1^{-1})\cdots (a_k b_k^{-1})\in D$. So $D$ is also a
group. As the order of $Y_\dt$ is $m$, $D$ is a subgroup of the
cyclic group $C_m$. Then $D=C_m$ as the order of $Y_\dt$ is $m$.

Let $\omega=\omega_m$ then in some suitable basis of
$\mathbb{C}^{n}$,
$$Y_\dt=diag(1,\cdots,1,\omega,\cdots,\omega,\cdots\cdots,\omega^{m-1},\cdots,\omega^{m-1}),$$
where for $j=0,1,\cdots,m-1,$ there are $t_j$ copies of $\omega^j$
on the diagonal with each $t_j>0$. Let $s=\frac{2\pi i}{m}$ and
$$A=diag(0,\cdots,0,s,\cdots,s,\cdots\cdots,(m-1)s,\cdots,(m-1)s),$$ where for $j=0,1,\cdots,m-1$ there are $t_j$ copies of $js$ on the
diagonal. Then $exp(A)=Y_\dt$.

Let $D=(d_{i,j})_{m\times m}$ where $d_{i,j}=w^{ij}$ for
$i,j=0,1,\cdots,m-1$. As $D$ is invertible, there are unique complex
numbers $c_0,c_1,\cdots,c_{m-1}$ satisfying
$$D\cdot (c_0,c_1,\cdots,c_{m-1})^{t}=(0,s,\cdots,(m-1)s)^{t}.$$
 Then
$\sum_{i=0}^{m-1}c_i Y_\dt^i=A$ and $exp(\sum_{i=0}^{m-1}c_i
Y_\dt^i)=Y_\dt.$ Then as $$\sum_{i=1}^{m-1}c_i\ {Y_\dt^i}\in
L_{[\dt]}\ \an \ \Lie Z(\Ker \dt)_0=\ad L_{[\dt]}$$ one has
$$\Ad_{Y_\dt}=exp(\ad(\sum_{i=1}^{m-1}c_i\ {Y_\dt^i}))\in Z(\Ker
\dt)_0.$$ \ep

Let $K$ be a finite maximal diagonalizable subgroup of
$G=\PGL(n,\mathbb{C})$. Recall that $K$ has a $W_G(K)$-invariant
anti-symmetric pairing $<,>$ and $(K,<,>)$ is a nonsingular
symplectic abelian group by Proposition \ref{b1}. Thus
$W_G(K)\st\Sp(K)$, where $\Sp(K)$ is the isometry group of
$(K,<,>)$.
\begin{theorem}\lb{d7}
 The Weyl group $W_G(K)$ equals $\Sp(K)$, and is generated by
the set of transvections $ s_\ga$ with $\ga\in\ddt(sl(n,\mb C),K)$.
\end{theorem}
\bp  By Lemma \ref{c1} one has $s_\ga\in W_G(K)$ for each
$\ga\in\ddt(sl(n,C),K)$. As $\ddt(sl(n,C),K)=\wh K\setminus \{0\}$,
$W_G(K)$ contains all $ s_\si$ with $\si$ a nonidentity element in
$K$. By Theorem \ref{c3} all such $ s_\si$ generate $\Sp(K)$, thus
$W_G(K)=\Sp(K)$.

\ep

\end{document}